\documentclass[10pt]{amsart}
\usepackage{amsfonts} 
\usepackage{graphicx,latexsym,bm,amsmath,amssymb,verbatim,multicol,lscape}
\theoremstyle{definition}

\theoremstyle{remark}

\numberwithin{equation}{section}

\begin{document}
\title{On deep holes of generalized Reed-Solomon codes}
\author{Shaofang Hong}
\address{Mathematical College, Sichuan University, Chengdu 610064, P.R. China}
\email{sfhong@scu.edu.cn, s-f.hong@tom.com, hongsf02@yahoo.com}
\author{Rongjun Wu}
\address{Mathematical College, Sichuan University, Chengdu 610064, P.R. China}
\email{eugen\_woo@163.com}
\thanks{Hong is the corresponding author and was supported partially
by National Science Foundation of China Grant \# 11371260 and by the
Ph.D. Programs Foundation of Ministry of Education of China Grant
\#20100181110073}

\date{\today}%
\begin{abstract}
Determining deep holes is an important topic in decoding
Reed-Solomon codes. In a previous paper \cite{[WH]}, we showed that
the received word $u$ is a deep hole of the standard Reed-Solomon
codes $[q-1, k]_q$ if its Lagrange interpolation polynomial is
the sum of monomial of degree $q-2$ and a polynomial of degree
at most $k-1$. In this paper, we extend this result by giving
a new class of deep holes of the generalized Reed-Solomon codes.
\end{abstract}

\maketitle

\section{\bf Introduction and the statement of the main result}

Let $\mathbb{F}_q$ be the finite field of $q$ elements with
characteristic $p$. Let $n$ and $k$ be positive integers.
Let $D=\{x_1, ..., x_n\}$ be a subset of $\mathbb{F}_q$,
which is called the {\it evaluation set}. The {\it generalized
Reed-Solomon code} $\mathcal {C}_q(D,k)$ of length $n$
and dimension $k$ over $\mathbb{F}_q$ is defined as follows:
$$
\mathcal {C}_q(D,k)=\{(f(x_1), ..., f(x_n))\in \mathbb{F}_q^n | f(x)
\in \mathbb{F}_q[x], {\rm deg}(f(x))\leq k-1\}.
$$
If $D=\mathbb{F}_q^*$, then it is called {\it standard Reed-Solomon code}.
If $D=\mathbb{F}_q$, then it is called {\it extended Reed-Solomon code}.
For any $[n,k]_q$ linear code $\mathcal {C}$, the
{\it minimum distance} $d(\mathcal {C})$ is defined by
$$
d(\mathcal {C}):= {\min}\{d(x,y)|x \in \mathcal{C}, y\in\mathcal{C}, x\neq y\},
$$
where $d(\cdot,\cdot)$ denotes the {\it Hamming distance} of two
words which is the number of different entries of them and
$w(\cdot)$ denotes the {\it Hamming weight} of a word which
is the number of its nonzero entries. Thus we have
$$
d(\mathcal {C})={\min}\{d(x,0)|0 \neq x \in \mathcal {C}\}
={\min}\{w(x)|0 \neq x \in \mathcal {C}\}.
$$
The {\it error distance} to code $\mathcal {C}$ of a received word
$u \in \mathbb{F}_q^n$ is defined by
$$d(u, \mathcal {C}):={\rm min}\{d(u, v)|v \in \mathcal {C}\}.$$
Clearly $d(u,\mathcal {C})=0$
if and only if $u \in \mathcal {C}$. The {\it covering radius}
$\rho(\mathcal {C})$ of code $\mathcal {C}$ is defined to be ${\rm
max}\{d(u,\mathcal {C})|u \in \mathbb{F}_p^n\}$. For the generalized
Reed-Solomon code $\mathcal {C}=\mathcal {C}_q(D,k)$, we have that
the minimum distance $d(\mathcal {C})=n-k+1$ and the covering radius
$\rho(\mathcal {C})=n-k$. The most important algorithmic problem in
coding theory is the maximum likelihood decoding (MLD): Given a
received word, find a word $v \in \mathcal {C}$ such that
$d(u,v)=d(u,\mathcal {C})$ [5]. Therefore, it is very
crucial to decide $d(u,\mathcal {C})$ for the word $u$. Sudan
[6] and Guruswami-Sudan [2] provided a polynomial
time list decoding algorithm for the decoding of $u$ when
$d(u,\mathcal {C})\leq n-\sqrt{nk}$. When the error distance
increases, the decoding becomes NP-complete for the generalized
Reed-Solomon codes [3].

When decoding the generalized Reed-Solomon code $\mathcal C$, for a
received word $u=(u_1, ..., u_n)\in \mathbb{F}_q^n$, we define the
{\it Lagrange interpolation polynomial} $u(x)$ of $u$ by
$$
u(x):=\sum_{i=1}^n u_i \prod_{\substack{j =1\\j \neq
i}}^{n}\frac{x-x_j} {x_i-x_j} \in \mathbb{F}_q[x],
$$
i.e., $u(x)$ is the unique polynomial of degree at most $n-1$ such
that $u(x_i)=u_i$ for $1 \leq i \leq n$. For $u \in \mathbb{F}_q^n$, we
define the degree of $u(x)$ to be the {\it degree} of $u$, i.e.,
${\rm deg}(u) ={\rm deg}(u(x))$. It is clear that $d(u,\mathcal
{C})=0$ if and only if ${\rm deg}(u) \leq k-1$. Evidently, we have
the following simple bounds.\\
\\
{\bf Lemma 1.1.} [4] {\it For $k \leq {\rm deg}(u)
\leq n-1$, we have the inequality
$$
n-{\rm deg}(u) \leq d(u,\mathcal {C}) \leq n-k = \rho.
$$}

Let $u \in \mathbb{F}_q^n$. If $d(u,\mathcal {C})=n-k,$ then the word
$u$ is called a {\it deep hole}. If ${\rm deg}(u)=k$, then the upper
bound is equal to the lower bound, and so $d(u,\mathcal {C})=n-k$
which implies that $u$ is a deep hole. This gives immediately
$(q-1)q^k$ deep holes. We call these deep holes {\it the trivial }
deep holes. It is an interesting open problem to determine all deep
holes. Cheng and Murray [1] showed that for the standard
Reed-Solomon code $[p-1, k]_p$ with $k < p^{1/4-\epsilon}$, the
received vector $(f(\alpha ))_{\alpha\in \mathbb{F}_p^*}$ cannot be a
deep hole if $f(x)$ is a polynomial of degree $k+d$ for $1\le
d<p^{3/13-\epsilon}$. Based on this result, they conjectured that
there is no other deep holes except the trivial ones mentioned
above. Li and Wan [5] used the method of character sums to
obtain a bound on the non-existence of deep holes for the extended
Reed-Solomon code $\mathcal {C}_q(\mathbb{F}_q,k)$. Wu and Hong [8]
found a counterexample to the Cheng-Murray conjecture
[1] about the standard Reed-Solomon codes.

Let $l$ be a positive integer. In this paper, we investigate
the deep holes of the generalized Reed-Solomon codes with
the evaluation set
$D:=\mathbb{F}_q \setminus\{a_1, ..., a_l\}$,
where $a_1, ..., a_l$ are any fixed $l$ distinct
elements of $\mathbb{F}_q$. Our method here
is different from that of [8].
Write $D=\{x_1, ..., x_{q-l}\}$
and for any $f(x) \in \mathbb{F}_q[x]$, let
$$f(D):=(f(x_1), ..., f(x_{q-l})).$$
Then we can rewrite the generalized Reed-Solomon code
$\mathcal {C}_q(D,k)$ with evaluation set $D$ as
$$
\mathcal {C}_q(D,k)=\{f(D)\in \mathbb{F}_q^{q-l} |
f(x) \in \mathbb{F}_q[x], {\rm deg}(f(x))\leq k - 1\}.
$$
Actually, by constructing some suitable auxiliary polynomials,
we find a new class of deep holes for the generalized
Reed-Solomon codes. That is, we have the following result.\\
\\
{\bf Theorem 1.2.} {\it Let $q \geq 4$ and $2 \leq k \leq q-l-1$.
For $1\le j\le l$, we define
\begin{equation} \label{1.1}
u_j(x):=\lambda_j(x-a_j)^{q-2}+r_j(x),
\end{equation}
where $\lambda_j \in \mathbb{F}^*_q$ and $r_j(x)\in \mathbb{F}_q[x]$
is a polynomial of degree at most $k-1$.
Then the received words $u_1(D), ..., u_l(D)$ are deep holes
of the generalized Reed-Solomon code $\mathcal {C}_q(D,k)$.}\\

The proof of Theorem 1.2 will be given in Section 2.

The materials presented here form part of the second author's
PhD thesis [7], which was finished on April 15, 2012.

\section{\bf Proof of Theorem 1.2}

Evidently, for any $a\in \mathbb{F}_q$, we have
$$\Big(\prod_{i=1}^{q-l}(a-x_i)\Big)\prod_{j=1}^{l}(a-a_j)=a^q-a=0,$$
and for any $a\in D$, we have $N(a)=0$, where
$$
N(x):=\prod_{i=1}^{q-l}(x-x_i).
$$
For $f(x) \in \mathbb{F}_q[x]$, by $\bar{f}(x)\in \mathbb{F}_q[x]$
we denote the reduction of $f(x) \mod N(x)$.
Therefore, for any $x_i \in D$, we have $f(x_i)=\bar{f}(x_i).$

First of all, we give a lemma about error distance.
In what follows, we let $G_k$ denote the set of all
the polynomials in $\mathbb{F}_q[x]$ of degree
at most $k-1$.\\
\\
{\bf Lemma 2.1.} \label{lm:uv} {\it Let $\#(D)=n$ and let
$u, v \in \mathbb{F}_q^n$ be two words.
If $u=\lambda v+f_{\leq k-1}(D)$, where $\lambda \in \mathbb{F}_q^*$
and $f_{\le k-1}(x)\in \mathbb{F}_q[x]$ is a polynomial of degree
at most $k-1$, then
$$
d(u, \mathcal{C}_q(D,k)) = d(v, \mathcal{C}_q(D,k)).
$$
Furthermore, $u$ is a deep hole of
$\mathcal{C}_q(D,k)$ if and only if $v$ is a deep hole of
$\mathcal{C}_q(D,k)$.}

\begin{proof} From the definition of error distance and
noting that $f_{\le k-1}(x)\in G_k$, we get immediately that
\begin{align*}
&d(u, \mathcal{C}_q(D,k))\\
=& \min_{g(x)\in G_k} \{d(u, g(D))\}\\
=& \min_{g(x)\in G_k} d(\lambda v + f_{\leq k-1}(D), g(D))\\
=& \min_{g(x)\in G_k} d(\lambda v + f_{\leq k-1}(D), g(D)+ f_{\leq k-1}(D))\\
=& \min_{g(x)\in G_k} d(\lambda v, g(D))\\
=& \min_{g(x)\in G_k} d(\lambda v, \lambda g(D)) \ \ \text{({\rm since}  $\lambda\ne 0$)}\\
=& \min_{g(x)\in G_k} d(v, g(D))\\
=& d(v, \mathcal{C}_q(D,k))
\end{align*}
as one desires. So Lemma 2.1 is proved.
\end{proof}

Now we are in the position to prove Theorem 1.2.\\
\\
{\it Proof of Theorem 1.2.} Let $f(x), g(x) \in \mathbb{F}_q[x]$.
One can deduce that
\begin{align}
\nonumber &d(f(D),g(D))\\
\nonumber =&\#\{x_i \in D \mid f(x_i) \ne g(x_i)\}\\
\nonumber =&\#\{x_i \in D \mid f(x_i) - g(x_i) \ne 0\}\\
=&\#(D) - \#\{x_i \in D \mid f(x_i) - g(x_i) = 0\}.\label{2.1}
\end{align}
Then by (\ref{2.1}), we infer that
\begin{align}
\nonumber   &d(f(D),\mathcal{C}_q(D,k))\\
\nonumber =&\min_{h(x)\in G_k}d(f(D),h(D))\\
\nonumber =&\min_{h(x)\in G_k}\{\#(D) - \#\{x_i \in D \mid f(x_i) - h(x_i) = 0\}\}\\
          =&q-l-\max_{h(x)\in G_k}\#\{x_i \in D \mid f(x_i) - h(x_i) = 0\}. \label{2.2}
\end{align}

For any integer $j$ with $1\le j\le l$, we let
$$f_j(x):=(x-a_j)^{q-2}\in \mathbb{F}_q[x].$$
For any $y \in D$, we have $y-a_j \neq 0$,
and so $f_j(y)=\frac{1}{y-a_j}$. We claim that
\begin{equation}\label{2.3}
\max_{h(x)\in G_k}\#\{y \in D \mid f_j(y) - h(y) = 0\}=k.
\end{equation}

In order to prove this claim, we pick $k$ distinct nonzero elements
$c_{j_1}, ..., c_{j_k}$ of $\mathbb{F}_q\setminus\{a_t-a_j\}_{t=1}^l$
(since $k\le q-l-1$). Now we introduce the auxiliary polynomial $g_j(x)$
as follows:
$$
g_j(x)=\frac{1}{x}\Big(1-\prod_{i=1}^{k}(1-c_{j_i}^{-1}x)\Big)\in \mathbb{F}_q[x].
$$
Then $\deg (g_j(x))=k-1$, and so $g_j(x)\in G_k$.
Since for any $y \in D$, we have
\begin{equation*}
\begin{split}
&f_j(y)-g_j(y-a_j)\\
=&\frac{1}{y-a_j}-g_j(y-a_j)\\
=&\frac{1}{y-a_j}(1-(y-a_j)g_j(y-a_j))\\
=&\frac{1}{y-a_j}\prod_{i=1}^{k}(1-c_{j_i}^{-1}(y-a_j)).
\end{split}
\end{equation*}
It then follows that $c_{j_1}+a_j, ..., c_{j_k}+a_j$
are the all roots of $f_j(x)-g_j(x-a_j)=0$ over
$\mathbb{F}_q$. Noticing that
$c_{j_1}, ..., c_{j_k} \in \mathbb{F}_q\setminus\{a_1-a_j, ..., a_l-a_j\}$,
we have $c_{j_1}+a_j, ..., c_{j_k}+a_j \in D$.
Also $D\subseteq \mathbb{F}_q$. Therefore
$c_{j_1}+a_j, ..., c_{j_k}+a_j$ are the
all roots of $f_j(x)-g_j(x-a_j)=0$ over $D$. Hence
\begin{equation}\label{2.4}
\#\{y \in D \mid f_j(y)-g_j(y-a_j)=0\}=k.
\end{equation}

On the other hand, for any $ h(x)\in G_k$, the equation
$1-(x-a_j)h(x)=0$ has at most $k$ roots over $\mathbb{F}_q$,
and so it has at most $k$ roots over $D$. But
$\frac{1}{y-a_j} \neq 0$ for any $y\in D$. Thus
\begin{equation*}
\begin{split}
&f_j(y)-h(y-a_j)\\
=&\frac{1}{y-a_j}-h(y-a_j)\\
=&\frac{1}{y-a_j}(1-(y-a_j)h(y-a_j)).
\end{split}
\end{equation*}
Hence for any $ h(x)\in G_k$, we have
\begin{equation*}
\#\{y \in D \mid f_j(y) - h(y) = 0\} \le k
\end{equation*}
which implies that
\begin{equation}\label{2.5}
\max_{h(x)\in G_k}\#\{y \in D \mid f_j(y)-h(y)=0\}\le k.
\end{equation}
From (\ref{2.4}) and (\ref{2.5}), we arrive at the desired
result (\ref{2.3}). The claim (\ref{2.3}) is proved.

Now from (\ref{2.2}) and (\ref{2.3}),  
we derive immediately that
$$d(f_j(D),\mathcal{C}_q(D,k)) = q-l-k.$$
In other words, $f_j(D)$ is a deep hole of the
generalized Reed-Solomon $\mathcal{C}_q(D,k)$.

Finally, from (\ref{1.1}) one can deduce that
\begin{equation}\label{2.6}
u_j(D)=\lambda_j f_j(D)+r_j(D).
\end{equation}
Since $\deg r_j(x)\le k-1$, it then follows
from (\ref{2.6}) and Lemma 2.1 that $u_j(D)$ is a deep
hole of $\mathcal{C}_q(D,k)$ as required.

This completes the proof of Theorem 1.2.
\hfill$\Box$\\

{\bf Acknowledgement}
The authors would like to thank the anonymous referee
for very careful reading of the manuscript and helpful
comments.

\end{document}